\documentclass[12pt,a4paper]{article}
\usepackage{exscale}
\usepackage{amsmath}
\usepackage{setspace}
\usepackage{amsfonts}
\usepackage{amssymb}
\usepackage{latexsym}

\newtheorem{lemma}{Lemma}[section]
\newtheorem{theorem}{Theorem}[section]
\newtheorem{corollary}{Corollary}[section]
\newtheorem{definicion}{Definition}[section]
\newtheorem{proposicion}{Proposition}[section]

\newtheorem{remark}{Remark}[section]

\begin{document}

\title{ Characterization of the minimal penalty of a convex risk measure
with applications to L\'{e}vy processes.}
\author{ {Daniel Hern\'{a}ndez--Hern\'{a}ndez\thanks{%
Centro de Investigaci\'{o}n en Matem\'{a}ticas, Apartado postal 402,
Guanajuato, Gto. 36000, M\'{e}xico. E-mail: dher@cimat.mx } \ \ Leonel P\'{e}%
rez-Hern\'{a}ndez\thanks{%
Departamento de Econom\'{\i}a y Finanzas, Universidad de Guanajuato, DCEA
Campus Guanajuato, C.P. 36250, Guanajuato, Gto. E-mail:
lperezhernandez@yahoo.com}}}
\maketitle

\begin{abstract}
The minimality of the penalty function associated with a convex risk
measure is analyzed in this paper. First, in a general static framework, we
provide necessary and sufficient conditions for a penalty function defined
in a convex and closed subset of the absolutely continuous measures with
respect to some reference measure $\mathbb{P}$ to be minimal on this set.
When the probability space supports a L\'{e}vy process, we establish results
that guarantee the minimality property of a penalty function described in
terms of the coefficients associated with the density processes. The set of
densities processes is described and the convergence of its quadratic
variation is analyzed.
\end{abstract}



\section{Introduction}

The definition of coherent risk measure was introduced by Artzner \textit{et
al.} in their fundamental works \cite{ADEH 1997}, \cite{ADEH 1999} for
finite probability spaces, giving an axiomatic characterization that was
extended later by Delbaen \cite{Delbaen 2002} to general probability spaces.
In the papers  mentioned above one of the fundamental axioms was the positive
homogeneity, and in further works it was removed, defining the concept of
convex risk measure introduced by F\"{o}llmer and Schied \cite{FoellSch 2002
a}, \cite{FoellSch 2002 b}, Frittelli and Rosazza Gianin \cite{FritRsza 2002}%
, \cite{FritRsza 2004} and Heath \cite{Heath 2000}.

This is a rich area that has received a lot of attention and much work has
been developed. There exists by now a well established theory in the static
and dynamic cases, but there are still many questions unanswered in the
static framework that need to be analyzed carefully. The one we focus on in
this paper is the characterization of the penalty functions that are minimal
for the corresponding static risk measure. Up to now, there are mainly two
ways to deal with minimal penalty functions, namely the definition or
the biduality relation. With the results presented in this paper we can
start with a penalty function, which essentially discriminate models within
a convex closed subset of absolutely continuous probability measures with
respect to (w.r.t.) the market measure, and then guarantee that it
corresponds to the minimal penalty of the corresponding convex risk measure
on this subset. This property is, as we will see, closely related with the
lower semicontinuity of the penalty function, and the complications to prove
this property depend on the structure of the probability space.

We first provide a general framework, within a measurable space with a
reference probability measure $\mathbb{P}$, and show necessary and
sufficient conditions for a penalty function defined in a convex and closed
subset of the absolutely continuous measures with respect to the reference
measure to be minimal within this subset. The characterization of the form of the
penalty functions that are minimal when the probability space supports a L%
\'{e}vy process is then studied. This requires to characterize the set of
absolutely continuous measures for this space, and it is done using results
that describe the density process for spaces which support semimartingales
with the weak predictable representation property. Roughly speaking, using
the weak representation property, every density process splits in two parts,
one is related with the continuous local martingale part of the
decomposition and the other with the corresponding discontinuous one. It is
shown some kind of continuity property for the quadratic variation of a
sequence of densities converging in \ $L^{1}$. From this characterization of
the densities,   a family of penalty functions is proposed, which turned out to
be minimal for the risk measures  generated by duality.

The paper is organized as follows. Section 2 contains the description of the
minimal penalty functions for a general probability space, providing
necessary and sufficient conditions, the last one rectricted to a subset of
equivalent probability measures. Section 3 reports the structure of the
densities for a probability space that supports a L\'{e}vy processes and the
convergence properties needed to prove the lower semicontinuity of the set
of penalty functions defined in Section 4. In this last section we show that
these penalty functions are minimal.

\section{Minimal penalty function of risk measures concentrated in $\mathcal{%
Q}_{\ll }\left( \mathbb{P}\right) $. \label{Sect Minimal Penalty Function of
CMR}}

\setcounter{equation}{0} Any penalty function $\psi $ induce a convex risk
measure $\rho $, which in turn has a representation by means of a minimal
penalty function $\psi _{\rho }^{\ast }.$ Starting with a penalty function $%
\psi $ concentrated in a convex and closed subset   of the set
of absolutely continuous probability measures with respect to some reference
measure $\mathbb{P}$, in this section we give necessary and sufficient conditions in order
to guarantee that $\psi $ is the minimal penalty within this set.  We begin recalling briefly some known
results from the theory of static risk measures, and then a characterization
for minimal penalties is presented.

\subsection{Preliminaries from static measures of risk \label%
{Subsect:_Preliminaries_SCRM}}

Let $X:\Omega \rightarrow \mathbb{R}$ be a mapping from a set $\Omega $ of
possible market scenarios, representing the discounted net worth of the
position. Uncertainty is represented by the measurable space $(\Omega,
\mathcal{F})$, and we denote by $\mathcal{X}$ the linear space of bounded
financial positions, including constant functions.

\begin{definicion}
\begin{enumerate}
\item[(i)] The function $\rho :\mathcal{X}\rightarrow \mathbb{R}$,
quantifying the risk of $X$, is a \textit{monetary risk measure} if it
satisfies the following properties:
\begin{equation}
\begin{array}{rl}
\text{Monotonicity:} & \text{If }X\leq Y\text{ then }\rho \left( X\right)
\geq \rho \left( Y\right) \ \forall X,Y\in \mathcal{X}.%
\end{array}
\label{Monotonicity}
\end{equation}%
$\smallskip \ $%
\begin{equation}
\begin{array}{rl}
\text{Translation Invariance:} & \rho \left( X+a\right) =\rho \left(
X\right) -a\ \forall a\in \mathbb{R}\ \forall X\in \mathcal{X}.%
\end{array}
\label{Translation Invariance}
\end{equation}

\item[(ii)] When this function satisfies also the convexity property
\begin{equation}
\begin{array}{rl}
& \rho \left( \lambda X+\left( 1-\lambda \right) Y\right) \leq \lambda \rho
\left( X\right) +\left( 1-\lambda \right) \rho \left( Y\right) \ \forall
\lambda \in \left[ 0,1\right] \ \forall X,Y\in \mathcal{X},%
\end{array}
\label{Convexity}
\end{equation}%
it is said that $\rho $ is a convex risk measure.

\item[(iii)] The function $\rho $ is called normalized if $\rho \left(
0\right) =0$, and sensitive, with respect to a measure $\mathbb{P}$, when
for each $X\in L_{+}^{\infty }\left( \mathbb{P}\right) $ with $\mathbb{P}%
\left[ X>0\right] >0$ we have that $\rho \left( -X\right) >\rho \left(
0\right) .$
\end{enumerate}
\end{definicion}

We say that a set function $%
\mathbb{Q}:\mathcal{F}\rightarrow \left[ 0,1\right] $ is a \textit{%
probability content} if it is finitely additive and $\mathbb{Q}\left( \Omega
\right) =1$. The set of \textit{probability contents} on this measurable
space is denoted by $\mathcal{Q}_{cont}$.
From the general theory of static convex risk measures \cite{FoellSch 2004}, we know that any map $%
\psi :\mathcal{Q}_{cont}\rightarrow \mathbb{R}\cup \{+\infty \},$ with $%
\inf\nolimits_{\mathbb{Q}\in \mathcal{Q}_{cont}}\psi (\mathbb{Q})\in \mathbb{%
R}$, induces a static convex measure of risk as a mapping $\rho :\mathfrak{M}%
_{b}\rightarrow \mathbb{R}$ given by
\begin{equation}
\rho (X):=\sup\nolimits_{\mathbb{Q}\in \mathcal{Q}_{cont}}\left\{ \mathbb{E}%
_{\mathbb{Q}}\left[ -X\right] -\psi (\mathbb{Q})\right\} .
\label{Static_CMR_induced_by_phi}
\end{equation}%
Here $\mathfrak{M}$ denotes the class of measurable functions and $\mathfrak{%
M}_{b}$ the subclass of bounded measurable functions. The function $\psi$
will be referred as a \textit{penalty function}. F\"{o}llmer and Schied \cite%
[Theorem 3.2]{FoellSch 2002 b} and Frittelli and Rosazza Gianin \cite[%
Corollary 7]{FritRsza 2002} proved that any convex risk measure is
essentially of this form.

More precisely, a convex risk measure  $\rho $ on the space  $\mathfrak{M}_{b}\left( \Omega ,\mathcal{F}\right) $ has the
representation
\begin{equation}
\rho (X)=\sup\limits_{\mathbb{Q}\in \mathcal{Q}_{cont}}\left\{ \mathbb{E}_{%
\mathbb{Q}}\left[ -X\right] -\psi _{\rho }^{\ast }\left( \mathbb{Q}\right)
\right\} ,  \label{Static_CMR_Robust_representation}
\end{equation}%
where
\begin{equation}
\psi _{\rho }^{\ast }\left( \mathbb{Q}\right) :=\sup\limits_{X\in \mathcal{A}%
\rho }\mathbb{E}_{\mathbb{Q}}\left[ -X\right] ,  \label{Def._minimal_penalty}
\end{equation}%
and $\mathcal{A}_{\rho }:=\left\{ X\in \mathfrak{M}_{b}:\rho (X)\leq
0\right\} $ is the \textit{acceptance set} of $\rho .$

The penalty $\psi _{\rho }^{\ast }$ is called the \textit{minimal penalty
function} associated to $\rho $ because, for any other penalty function $%
\psi $ fulfilling $\left( \ref{Static_CMR_Robust_representation}\right) $, $%
\psi \left( \mathbb{Q}\right) \geq \psi _{\rho }^{\ast }\left( \mathbb{Q}%
\right) $, for all $\mathbb{Q}\in \mathcal{Q}_{cont}.$ Furthermore, for the
minimal penalty function, the next biduality relation is satisfied
\begin{equation}
\psi _{\rho }^{\ast }\left( \mathbb{Q}\right) =\sup_{X\in \mathfrak{M}%
_{b}\left( \Omega ,\mathcal{F}\right) }\left\{ \mathbb{E}_{\mathbb{Q}}\left[
-X\right] -\rho \left( X\right) \right\} ,\quad \forall \mathbb{Q\in }%
\mathcal{Q}_{cont}.  \label{static convex rsk msr biduality}
\end{equation}

Let $\mathcal{Q}\left( \Omega ,\mathcal{F}\right) $ be the family of
probability measures on the measurable space $\left( \Omega ,\mathcal{F}%
\right) .$ Among the measures of risk, the class of them that are
concentrated on the set of probability measures $\mathcal{Q\subset Q}_{cont}$
are of special interest. Recall that a function $I:E\subset \mathbb{R}%
^{\Omega }\rightarrow \mathbb{R}$ is \textit{sequentially continuous from
below (above)} when $\left\{ X_{n}\right\} _{n\in \mathbb{N}}\uparrow
X\Rightarrow \lim_{n\rightarrow \infty }I\left( X_{n}\right) =I\left(
X\right) $ ( respectively $\left\{ X_{n}\right\} _{n\in \mathbb{N}%
}\downarrow X\Rightarrow \lim_{n\rightarrow \infty }I\left( X_{n}\right)
=I\left( X\right) $). F\"{o}llmer and Schied \cite{FoellSch 2004} proved
that any sequentially continuous from below convex measure of risk is
concentrated on the set $\mathcal{Q}$. Later, Kr\"{a}tschmer \cite[Prop. 3
p. 601]{Kraetschmer 2005} established that the sequential continuity from
below is not only a sufficient but also a necessary condition in order to
have a representation, by means of the minimal penalty function in terms of
probability measures.

We denote by $\mathcal{Q}_{\ll }(\mathbb{P})$ the subclass of absolutely
continuous probability measure with respect to $\mathbb{P}$ and by $\mathcal{%
Q}_{\approx }\left( \mathbb{P}\right) $ the subclass of equivalent
probability measure. Of course, $\mathcal{Q}_{\approx }\left( \mathbb{P}%
\right) \subset \mathcal{Q}_{\ll }(\mathbb{P})\subset \mathcal{Q}\left(
\Omega ,\mathcal{F}\right) $.

\begin{remark}\label{Remarkpsi(Q)=+oo_for_Q_not_<<} 
When a convex risk measures in $%
\mathcal{X}:=L^{\infty }\left( \mathbb{P}\right) $ satisfies  the
property
\begin{equation}
\rho \left( X\right) =\rho \left( Y\right) \text{ if }X=Y\ \mathbb{P}\text{%
-a.s.}  \label{rho(X)=rho(Y)_for_X=Y}
\end{equation}%
and is represented by a penalty function $\psi $ as in $\left( \ref%
{Static_CMR_induced_by_phi}\right) $, we have that
\begin{equation}
\mathbb{Q}\in \mathcal{Q}_{cont}\setminus \mathcal{Q}_{cont}^{\ll
}\Longrightarrow \psi \left( \mathbb{Q}\right) =+\infty ,
\label{psi(Q)=+oo_for_Q_not_<<}
\end{equation}%
where $\mathcal{Q}_{cont}^{\ll }$ is the set of contents absolutely
continuous with respect to $\mathbb{P}$; see \cite[Lemma 4.30 p. 172]%
{FoellSch 2004}.
\end{remark}

\subsection{Minimal penalty functions \label%
{Subsect:_Minimal_penalty_functions}}

The minimality property of the penalty function turns out to be quite
relevant, and it is a desirable property that is not easy to prove in
general. For instance, in the study of robust portfolio optimization
problems (see, for example, Schied \cite{Schd 2007} and Hern\'{a}ndez-Hern%
\'{a}ndez and P\'{e}rez-Hern\'{a}ndez \cite{PerHer}), using techniques of
duality, the minimality property is a necessary condition in order to have a
well posed dual problem. More recently, the dual representations of dynamic
risk measures were analyzed by Barrieu and El Karoui \cite{BaElKa2009},
while the connection with BSDEs and $g-$expectations have been studied by
Delbaen \textit{et. al.} \cite{DelPenRz}. The minimality of the penalty
function also plays a crucial role in the characterization of the time
consistency property for dynamic risk measures (see Bion-Nadal \cite%
{BionNa2008}, \cite{BionNa2009}).

In the next sections we will show some of the difficulties that appear to prove the minimality of the penalty function when
the probability space $(\Omega, \mathcal{F},\mathbb{P})$ supports a L\'evy
process. However, to establish the results of this section we only need to fix a
probability space $(\Omega, \mathcal{F}, \mathbb{P})$.

When we deal with a set of absolutely continuous probability measures $%
\mathcal{K}\subset \mathcal{Q}_{\ll }(\mathbb{P})$ it is necessary to make
reference to some topological concepts, meaning that we are considering the
corresponding set of densities and the strong topology in $L^{1}\left(
\mathbb{P}\right) .$ Recall that within a locally convex space, a convex set $%
\mathcal{K}$ is weakly closed if and only if $\mathcal{K}$ is closed in the
original topology \cite[Thm A.59]{FoellSch 2004}.

\begin{lemma}
\label{static minimal penalty funct. in Q(<<) <=>} Let $\psi :\mathcal{K}%
\subset \mathcal{Q}_{\ll }(\mathbb{P})\rightarrow \mathbb{R}\cup \{+\infty
\} $ be a function with $\inf\nolimits_{\mathbb{Q}\in \mathcal{K}}\psi (%
\mathbb{Q})\in \mathbb{R},$ and define the extension $\psi (\mathbb{Q}%
):=\infty $ for each $\mathbb{Q}\in \mathcal{Q}_{cont}\setminus \mathcal{K}$%
, with $\mathcal{K}$ a convex closed set. Also, define the function $\Psi $,
with domain in $L^{1}$, as
\begin{equation*}
\Psi \left( D\right) :=\left\{
\begin{array}{rl}
\psi \left( \mathbb{Q}\right) & \text{if }D=d\mathbb{Q}/d\mathbb{P}\text{
for }\mathbb{Q}\in \mathcal{K} \\
\infty & \text{otherwise.}%
\end{array}%
\right.
\end{equation*}%
Then, for the convex measure of risk $\rho (X):=\sup\limits_{\mathbb{Q}\in
\mathcal{Q}_{cont}}\left\{ \mathbb{E}_{\mathbb{Q}}\left[ -X\right] -\psi
\left( \mathbb{Q}\right) \right\} $ associated with $\psi $ the following assertions hold:

$\left( a\right) $ If $\rho $ has as minimal penalty $\psi _{\rho }^{\ast }$
the function $\psi $ (i.e. $\psi $ $=\psi _{\rho }^{\ast }$ ), then $\Psi $
is a proper convex function and lower semicontinuous w.r.t. the (strong) $%
L^{1}$-topology or equivalently w.r.t. the weak topology $\sigma \left(
L^{1},L^{\infty }\right) $.\

$\left( b\right) $ If $\Psi $ is lower semicontinuous w.r.t. the (strong) $%
L^{1}$-topology or equivalently w.r.t. the weak topology $\sigma \left(
L^{1},L^{\infty }\right) ,$ then
\begin{equation}
\psi \mathbf{1}_{\mathcal{Q}_{\ll }(\mathbb{P})}=\psi _{\rho }^{\ast }%
\mathbf{1}_{\mathcal{Q}_{\ll }(\mathbb{P})}.
\label{PSI_l.s.c=>psi*=psi_on_Q<<}
\end{equation}
\end{lemma}

\noindent \textit{Proof:} $\left( a\right) $ Recall that $\sigma \left(
L^{1},L^{\infty }\right) $ is the coarsest topology on $L^{1}\left( \mathbb{P%
}\right) $ under which every linear operator is continuous,  and hence
$\Psi _{0}^{X}\left( Z\right) :=\mathbb{E}_{\mathbb{P}}\left[ Z\left(
-X\right) \right] $, with $Z\in L^1$, is a continuous function for each $X\in
\mathfrak{M}_{b}\left( \Omega ,\mathcal{F}\right) $ fixed. For $\delta \left(
\mathcal{K}\right) :=\left\{ Z:Z=d\mathbb{Q}/d\mathbb{P}\text{ with }\mathbb{%
Q}\in \mathcal{K}\right\} $ we have that%
\begin{equation*}
\Psi _{1}^{X}\left( Z\right) :=\Psi _{0}^{X}\left( Z\right) \mathbf{1}%
_{\delta \left( \mathcal{K}\right) }\left( Z\right) +\infty \times \mathbf{1}%
_{L^{1}\setminus \delta \left( \mathcal{K}\right) }\left( Z\right)
\end{equation*}%
is clearly lower semicontinuous on $\delta \left( \mathcal{K}\right) .$ For $%
Z^{\prime }\in L^{1}\left( \mathbb{P}\right) \setminus \delta \left(
\mathcal{K}\right) $ arbitrary fixed we have from Hahn-Banach's Theorem that
there is a continuous lineal functional $l\left( Z\right) $ with $l\left(
Z^{\prime }\right) <\inf_{Z\in \delta \left( \mathcal{K}\right) }l\left(
Z\right) $. Taking $\varepsilon :=\frac{1}{2}\left\{ \inf_{Z\in \delta
\left( \mathcal{K}\right) }l\left( Z\right) -l\left( Z^{\prime }\right)
\right\} $ we have that the weak open ball $B\left( Z^{\prime },\varepsilon
\right) :=\left\{ Z\in L^{1}\left( \mathbb{P}\right) :\left\vert l\left(
Z^{\prime }\right) -l\left( Z\right) \right\vert <\varepsilon \right\} $
satisfies $B\left( Z^{\prime },\varepsilon \right) \cap \delta \left( \mathcal{%
K}\right) =\varnothing .$ Therefore, $\Psi _{1}^{X}\left( Z\right) $ is weak
lower semicontinuous on $L^{1}\left( \mathbb{P}\right) ,$ as well as $\Psi
_{2}^{X}\left( Z\right) :=\Psi _{1}^{X}\left( Z\right) -\rho \left( X\right)
.$ If $$\psi \left( \mathbb{Q}\right) =\psi _{\rho }^{\ast }\left( \mathbb{Q}%
\right) =\sup_{X\in \mathfrak{M}_{b}\left( \Omega ,\mathcal{F}\right)
}\left\{ \int Z\left( -X\right) d\mathbb{P}-\rho \left( X\right) \right\}, $$
where $Z:=d\mathbb{Q}/d\mathbb{P},$ we have that $\Psi \left( Z\right)
=\sup_{X\in \mathfrak{M}_{b}\left( \Omega ,\mathcal{F}\right) }\left\{ \Psi
_{2}^{X}\left( Z\right) \right\} $ is the supremum of a family of convex
lower semicontinuous functions with respect to the topology $\sigma \left(
L^{1},L^{\infty }\right) $, and $\Psi \left( Z\right) $ preserves both
properties.

$\left( b\right) $ For the Fenchel - Legendre transform (conjugate function)
$\Psi ^{\ast }:\ L^{\infty }\left( \mathbb{P}\right) \longrightarrow \mathbb{%
R}$  for each $U\in L^{\infty }\left( \mathbb{P}\right) $%
\begin{equation*}
\Psi ^{\ast }\left( U\right) =\sup\limits_{Z\in \delta \left( \mathcal{K}%
\right) }\left\{ \int ZUd\mathbb{P-}\Psi \left( Z\right) \right\}
=\sup\limits_{\mathbb{Q}\in \mathcal{Q}_{cont}}\left\{ \mathbb{E}_{\mathbb{Q}%
}\left[ U\right] \mathbb{-\psi }\left( \mathbb{Q}\right) \right\} \equiv
\rho \left( -U\right) .
\end{equation*}
From the lower semicontinuity of $\Psi $ w.r.t. the weak topology
$\sigma \left( L^{1},L^{\infty }\right) $  that $\Psi =\Psi ^{\ast \ast }$.
Considering the weak$^{\ast }$-topology $\sigma \left( L^{\infty }\left(
\mathbb{P}\right) ,L^{1}\left( \mathbb{P}\right) \right) $  for $Z=d%
\mathbb{Q}/d\mathbb{P}$ we have that
\begin{equation*}
\psi \left( \mathbb{Q}\right) =\Psi \left( Z\right) =\Psi ^{\ast \ast
}\left( Z\right) =\sup\limits_{U\in L^{\infty }\left( \mathbb{P}\right)
}\left\{ \int Z\left( -U\right) d\mathbb{P-}\Psi ^{\ast }\left( -U\right)
\right\} =\psi _{\rho }^{\ast }\left( \mathbb{Q}\right) .
\end{equation*}%
\hfill $\Box $

\begin{remark}
\begin{enumerate}
\item As pointed out in Remark \ref{Remarkpsi(Q)=+oo_for_Q_not_<<}, we have that $$\mathbb{Q}\in \mathcal{Q}%
_{cont}\setminus \mathcal{Q}_{cont}^{\ll }\Longrightarrow \psi _{\rho
}^{\ast }\left( \mathbb{Q}\right) =+\infty =\psi \left( \mathbb{Q}\right). $$
Therefore, under the conditions of Lemma \ref{static minimal penalty funct. in
Q(<<) <=>} $\left( b\right) $ the penalty function $\psi $ might differ from $\psi _{\rho
}^{\ast }$
on $\mathcal{Q}_{cont}^{\ll }\setminus \mathcal{Q}_{\ll }.$
For instance, the penalty function defined as $\psi \left( \mathbb{Q}\right) :=\infty \times \mathbf{1%
}_{\mathcal{Q}_{cont}\setminus  \mathcal{Q}_{\ll }(\mathbb{P})}\left( \mathbb{%
Q}\right) $ leads to the worst case risk measure $\rho (X):=\sup\nolimits_{%
\mathbb{Q}\in \mathcal{Q}_{\ll }(\mathbb{P})}\mathbb{E}_{\mathbb{Q}}\left[ -X%
\right] $, which has as minimal penalty the function $$\psi _{\rho }^{\ast
}\left( \mathbb{Q}\right) =\infty \times \mathbf{1}_{\mathcal{Q}%
_{cont}\setminus \mathcal{Q}_{cont}^{\ll }}\left( \mathbb{Q}\right). $$

\item  Note that the total variation distance $d_{TV}\left( \mathbb{Q}^{1},\mathbb{Q%
}^{2}\right) :=\sup_{A\in \mathcal{F}}\left\vert \mathbb{Q}^{1}\left[ A%
\right] -\mathbb{Q}^{2}\left[ A\right] \right\vert $, with $\mathbb{Q}^{1},\;\mathbb{Q}^{2}\in \mathcal{Q}_{\ll }$, fulfills that $%
d_{TV}\left( \mathbb{Q}^{1},\mathbb{Q}^{2}\right) \leq \left\Vert d\mathbb{Q}%
^{1}/d\mathbb{P}-\mathbb{Q}^{2}/d\mathbb{P}\right\Vert _{L^{1}}$. Therefore,
the minimal penalty function is lower semicontinuous in the total variation
topology; see Remark 4.16 (b) p. 163 in \cite{FoellSch 2004}.
\end{enumerate}
\end{remark}

\section{Preliminaries from stochastic calculus\label{Sect. Preliminaries}}

\setcounter{equation}{0} Within a probability space which supports a
semimartingale with the weak predictable representation property, there is a
representation of the density processes of the absolutely continuous
probability measures by means of two coefficients. Roughly speaking, this
means that the \textquotedblleft dimension\textquotedblright\ of the linear
space of local martingales is two. Throughout these coefficients we can
represent every local martingale as a combination of two components, namely
as an stochastic integral with respect to the continuous part of the
semimartingale and an integral with respect to its compensated jump measure.
This is of course the case for local martingales, and with more reason this
observation about the dimensionality holds for the martingales associated
with the corresponding densities processes. In this section we review those
concepts of stochastic calculus that are relevant to understand this
representation properties, and prove some kind of continuity property for
the quadratic variation of a sequence of uniformly integrable martingales
converging in \ $L^{1}$. This result is one of the contributions of this
paper.

\subsection{Fundamentals of L\'{e}vy and semimartingales processes \label%
{Subsect:_Fundamentals_Levy_and_Semimartingales}}

Let $\left( \Omega ,\mathcal{F},\mathbb{P}\right) $ be a probability space.
We say that $L:=\left\{ L_{t}\right\} _{t\in \mathbb{R}_{+}}$ is a L\'{e}vy
process for this probability space if it is an adapted c\`{a}dl\`{a}g
process with independent stationary increments starting at zero. The
filtration considered is $\mathbb{F}:=\left\{ \mathcal{F}_{t}^{\mathbb{P}%
}\left( L\right) \right\} _{t\in \mathbb{R}_{+}}$, the completion of its
natural filtration, i.e. $\mathcal{F}_{t}^{\mathbb{P}}\left( L\right)
:=\sigma \left\{ L_{s}:s\leq t\right\} \vee \mathcal{N}$ where $\mathcal{N}$
is the $\sigma $-algebra generated by all $\mathbb{P}$-null sets. The jump
measure of $L$ is denoted by $\mu :\Omega \times \left( \mathcal{B}\left(
\mathbb{R}_{+}\right) \otimes \mathcal{B}\left( \mathbb{R}_{0}\right)
\right) \rightarrow \mathbb{N}$ where $\mathbb{R}_{0}:=\mathbb{R}\setminus
\left\{ 0\right\} $. The dual predictable projection of this measure, also
known as its L\'{e}vy system, satisfies the relation $\mu ^{\mathcal{P}%
}\left( dt,dx\right) =dt\times \nu \left( dx\right) $, where $\nu \left(
\cdot \right) :=\mathbb{E}\left[ \mu \left( \left[ 0,1\right] \times \cdot
\right) \right] $ is the intensity or L\'{e}vy measure of $L.$

The L\'{e}vy-It\^{o} decomposition of $L$ is given by
\begin{equation}
L_{t}=bt+W_{t}+\int\limits_{\left[ 0,t\right] \times \left\{ 0<\left\vert
x\right\vert \leq 1\right\} }xd\left\{ \mu -\mu ^{\mathcal{P}}\right\}
+\int\limits_{\left[ 0,t\right] \times \left\{ \left\vert x\right\vert
>1\right\} }x\mu \left( ds,dx\right) .  \label{Levy-Ito_decomposition}
\end{equation}%
It implies that $L^{c}=W$ is the Wiener process, and hence $\left[ L^{c}%
\right] _{t}=t$, where $\left( \cdot \right) ^{c}$ and $\left[ \,\cdot \,%
\right] $ denote the continuous martingale part and the process of quadratic
variation of any semimartingale, respectively. For the predictable quadratic
variation we use the notation $\left\langle \,\cdot \,\right\rangle $.

Denote by $\mathcal{V}$ the set of c\`{a}dl\`{a}g, adapted processes with
finite variation, and let $\mathcal{V}^{+}\subset \mathcal{V}$ be the subset
of non-decreasing processes in $\mathcal{V}$ starting at zero. Let $\mathcal{%
A}\subset \mathcal{V}$ be the class of processes with integrable variation,
i.e. $A\in \mathcal{A}$ if and only if $\bigvee_{0}^{\infty }A\in
L^{1}\left( \mathbb{P}\right) $, where $\bigvee_{0}^{t}A$ denotes the
variation of $A$ over the finite interval $\left[ 0,t\right] $.
%\end{document}
The subset $\mathcal{A}^{+}=\mathcal{A\cap V}^{+}$ represents those
processes which are also increasing i.e. with non-negative right-continuous
increasing trajectories. Furthermore, $\mathcal{A}_{loc}$ (resp. $\mathcal{A}%
_{loc}^{+}$) is the collection of adapted processes with locally integrable
variation (resp. adapted locally integrable increasing processes). For a c%
\`{a}dl\`{a}g process $X$ we denote by $X_{-}:=\left( X_{t-}\right) $ the
left hand limit process, where $X_{0-}:=X_{0}$ by convention, and by $%
\bigtriangleup X=\left( \bigtriangleup X_{t}\right) $ the jump process $%
\bigtriangleup X_{t}:=X_{t}-X_{t-}$.

Given an adapted c\`{a}dl\`{a}g semimartingale $U$, the jump measure and its
dual predictable projection (or compensator) are denoted by $\mu _{U}\left( %
\left[ 0,t\right] \times A\right) :=\sum_{s\leq t}\mathbf{1}_{A}\left(
\triangle U_{s}\right) $ and $\mu _{U}^{\mathcal{P}}$, respectively.
Further, we denote by $\mathcal{P}\subset \mathcal{F}\otimes \mathcal{B}%
\left( \mathbb{R}_{+}\right) $ the predictable $\sigma $-algebra and by $%
\widetilde{\mathcal{P}}:=\mathcal{P}\otimes \mathcal{B}\left( \mathbb{R}%
_{0}\right) .$ With some abuse of notation, we write $\theta _{1}\in
\widetilde{\mathcal{P}}$ when the function $\theta _{1}:$ $\Omega \times
\mathbb{R}_{+}\times \mathbb{R}_{0}\rightarrow \mathbb{R}$ is $\widetilde{%
\mathcal{P}}$-measurable and $\theta \in \mathcal{P}$ for predictable
processes.

Let
\begin{equation}
\begin{array}{clc}
\mathcal{L}\left( U^{c}\right) := & \left\{ \theta \in \mathcal{P}:\exists
\left\{ \tau _{n}\right\} _{n\in \mathbb{N}}\text{ sequence of stopping
times with }\tau _{n}\uparrow \infty \right. &  \\
& \left. \text{and }\mathbb{E}\left[ \int\limits_{0}^{\tau _{n}}\theta ^{2}d%
\left[ U^{c}\right] \right] <\infty \ \forall n\in \mathbb{N}\right\} &
\end{array}
\label{Def._L(U)}
\end{equation}%
be the class of predictable processes $\theta \in \mathcal{P}$ integrable
with respect to $U^{c}$ in the sense of local martingale, and by
\begin{equation*}
\Lambda \left( U^{c}\right) :=\left\{ \int \theta _{0}dU^{c}:\theta _{0}\in
\mathcal{L}\left( U^{c}\right) \right\}
\end{equation*}%
the linear space of processes which admits a representation as the
stochastic integral with respect to $U^{c}$. For an integer valued random
measure $\mu ^{\prime }$ we denote by $\mathcal{G}\left( \mu ^{\prime
}\right) $ the class of $\widetilde{\mathcal{P}}$-measurable processes $%
\theta _{1}:$ $\Omega \times \mathbb{R}_{+}\times \mathbb{R}_{0}\rightarrow
\mathbb{R}$ satisfying the following conditions:
\begin{equation*}
\begin{array}{cl}
\left( i\right) & \theta _{1}\in \widetilde{\mathcal{P}}, \\
\left( ii\right) & \int\limits_{\mathbb{R}_{0}}\left\vert \theta _{1}\left(
t,x\right) \right\vert \left( \mu ^{\prime }\right) ^{\mathcal{P}}\left(
\left\{ t\right\} ,dx\right) <\infty \ \forall t>0, \\
\left( iii\right) & \text{The process } \\
& \left\{ \sqrt{\sum\limits_{s\leq t}\left\{ \int\limits_{\mathbb{R}%
_{0}}\theta _{1}\left( s,x\right) \mu ^{\prime }\left( \left\{ s\right\}
,dx\right) -\int\limits_{\mathbb{R}_{0}}\theta _{1}\left( s,x\right) \left(
\mu ^{\prime }\right) ^{\mathcal{P}}\left( \left\{ s\right\} ,dx\right)
\right\} ^{2}}\right\} _{t\in \mathbb{R}_{+}}\in \mathcal{A}_{loc}^{+}.%
\end{array}%
\end{equation*}%
The set $\mathcal{G}\left( \mu ^{\prime }\right) $ represents the domain of
the functional $\theta _{1}\rightarrow \int \theta _{1}d\left( \mu ^{\prime
}-\left( \mu ^{\prime }\right) ^{\mathcal{P}}\right) ,$ which assign to $%
\theta _{1}$ the unique purely discontinuous local martingale $M$ with
$$
\bigtriangleup M_{t}=\int\limits_{\mathbb{R}_{0}}\theta _{1}\left(
t,x\right) \mu ^{\prime }\left( \left\{ t\right\} ,dx\right) -\int\limits_{%
\mathbb{R}_{0}}\theta _{1}\left( t,x\right) \left( \mu ^{\prime }\right) ^{%
\mathcal{P}}\left( \left\{ t\right\} ,dx\right) .$$

We use the notation $\int
\theta _{1}d\left( \mu ^{\prime }-\left( \mu ^{\prime }\right) ^{\mathcal{P}%
}\right) $ to write the value of this functional in $\theta _{1}$. It is
important to point out that this functional is not, in general, the integral
with respect to the difference of two measures. For a detailed exposition on
these topics see He, Wang and Yan \cite{HeWanYan} or Jacod and Shiryaev \cite%
{Jcd&Shry 2003}, which are our basic references.

In particular, for the L\'{e}vy process $L$ with jump measure $\mu $,
\begin{equation}
\mathcal{G}\left( \mu \right) \equiv \left\{ \theta _{1}\in \widetilde{%
\mathcal{P}}:\left\{ \sqrt{\sum\limits_{s\leq t}\left\{ \theta _{1}\left(
s,\triangle L_{s}\right) \right\} ^{2}\mathbf{1}_{\mathbb{R}_{0}}\left(
\triangle L_{s}\right) }\right\} _{t\in \mathbb{R}_{+}}\in \mathcal{A}%
_{loc}^{+}\right\} ,  \label{G(miu) Definition}
\end{equation}%
since $\mu ^{\mathcal{P}}\left( \left\{ t\right\} \times A\right) =0$, for
any Borel set $A$ of $\mathbb{R}_{0}$.

We say that the semimartingale $U$ has the \textit{weak property of
predictable representation} when
\begin{equation}
\mathcal{M}_{loc,0}=\Lambda \left( U^{c}\right) +\left\{ \int \theta
_{1}d\left( \mu _{U}-\mu _{U}^{\mathcal{P}}\right) :\theta _{1}\in \mathcal{G%
}\left( \mu _{U}\right) \right\} ,\   \label{Def_weak_predictable_repres.}
\end{equation}%
where the previous sum is the linear sum of the vector spaces, and $\mathcal{%
M}_{loc,0}$ is the linear space of local martingales starting at zero.

Let $\mathcal{M}$ and $\mathcal{M}_{\infty }$ denote the class of c\`{a}dl%
\`{a}g and c\`{a}dl\`{a}g uniformly integrable martingale respectively. The
following lemma is interesting by itself to understand the continuity
properties of the quadratic variation for a given convergent sequence of
uniformly integrable martingale . It will play a central role in the proof
of the lower semicontinuity of the penalization function introduced in
section \ref{Sect Penalty Function for densities}. Observe that the
assertion of this lemma is valid in a general filtered probability space and
not only for the completed natural filtration of the L\'{e}vy process
introduced above.

\begin{lemma}
\label{E[|Mn-M|]->0=>[Mn-M](oo)->0_in_P}For $\left\{ M^{\left( n\right)
}\right\} _{n\in \mathbb{N}}\subset \mathcal{M}_{\infty }$ and $M\in
\mathcal{M}_{\infty }$ the following implication holds
\begin{equation*}
M_{\infty }^{\left( n\right) }\overset{L^{1}}{\underset{n\rightarrow \infty }%
{\longrightarrow }}M_{\infty }\Longrightarrow \left[ M^{\left( n\right) }-M%
\right] _{\infty }\overset{\mathbb{P}}{\longrightarrow }0.
\end{equation*}%
Moreover,%
\begin{equation*}
M_{\infty }^{\left( n\right) }\overset{L^{1}}{\underset{n\rightarrow \infty }%
{\longrightarrow }}M_{\infty }\Longrightarrow \left[ M^{\left( n\right) }-M%
\right] _{t}\overset{\mathbb{P}}{\underset{n\rightarrow \infty }{%
\longrightarrow }}0\;\; \forall t.
\end{equation*}
\end{lemma}

\noindent \textit{Proof.} From the $L^{1}$ convergence of $M_{\infty
}^{\left( n\right) }$ to $M_{\infty }$, we have that $\{M_{\infty }^{\left(
n\right) }\}_{n\in \mathbb{N}}\cup \left\{ M_{\infty }\right\} $ is
uniformly integrable, which is equivalent to the existence of a convex and
increasing function $G:[0,+\infty )\rightarrow \lbrack 0,+\infty )$ such
that
\begin{equation*}
\left( i\right) \quad \lim_{x\rightarrow \infty }\frac{G\left( x\right) }{x}%
=\infty ,
\end{equation*}%
and
\begin{equation*}
\left( ii\right) \quad \sup_{n\in \mathbb{N}}\mathbb{E}\left[ G\left(
\left\vert M_{\infty }^{\left( n\right) }\right\vert \right) \right] \vee
\mathbb{E}\left[ G\left( \left\vert M_{\infty }\right\vert \right) \right]
<\infty .
\end{equation*}%
Now, define the stopping times
\begin{equation*}
\tau _{k}^{n}:=\inf \left\{ u>0:\sup_{t\leq u}\left\vert M_{t}^{\left(
n\right) }-M_{t}\right\vert \geq k\right\} .
\end{equation*}%
Observe that the estimation $\sup_{n\in \mathbb{N}}\mathbb{E}\left[ G\left(
\left\vert M_{\tau _{k}^{n}}^{\left( n\right) }\right\vert \right) \right]
\leq \sup_{n\in \mathbb{N}}\mathbb{E}\left[ G\left( \left\vert M_{\infty
}^{\left( n\right) }\right\vert \right) \right] $ implies the uniformly
integrability of $\left\{ M_{\tau _{k}^{n}}^{\left( n\right) }\right\}
_{n\in \mathbb{N}}$ for each $k$ fixed. Since any uniformly integrable c\`{a}%
dl\`{a}g martingale is of class $\mathcal{D}$, follows the uniform
integrability of $\left\{ M_{\tau _{k}^{n}}\right\} _{n\in \mathbb{N}}$ for
all $k\in \mathbb{N}$, and hence $\left\{ \sup\nolimits_{t\leq \tau
_{k}^{n}}\left\vert M_{t}^{\left( n\right) }-M_{t}\right\vert \right\}
_{n\in \mathbb{N}}$ is uniformly integrable. This and the maximal inequality
for supermartingales
\begin{eqnarray*}
\mathbb{P}\left[ \sup_{t\in \mathbb{R}_{+}}\left\vert M_{t}^{\left( n\right)
}-M_{t}\right\vert \geq \varepsilon \right] &\leq &\frac{1}{\varepsilon }%
\left\{ \sup_{t\in \mathbb{R}_{+}}\mathbb{E}\left[ \left\vert M_{t}^{\left(
n\right) }-M_{t}\right\vert \right] \right\} \\
&\leq &\frac{1}{\varepsilon }\mathbb{E}\left[ \left\vert M_{\infty }^{\left(
n\right) }-M_{\infty }\right\vert \right] \longrightarrow 0,
\end{eqnarray*}%
yields the convergence of $\left\{ \sup\nolimits_{t\leq \tau
_{k}^{n}}\left\vert M_{t}^{\left( n\right) }-M_{t}\right\vert \right\}
_{n\in \mathbb{N}}$ in $L^{1}$ to $0$. The second Davis' inequality \cite[Thm. 10.28]{HeWanYan} guarantees
that, for some constant $C$,
\begin{equation*}
\mathbb{E}\left[ \sqrt{\left[ M^{\left( n\right) }-M\right] _{\tau _{k}^{n}}}%
\right] \leq C\mathbb{E}\left[ \sup\limits_{t\leq \tau _{k}^{n}}\left\vert
M_{t}^{\left( n\right) }-M_{t}\right\vert \right] \underset{n\rightarrow
\infty }{\longrightarrow }0\quad \forall k\in \mathbb{N},
\end{equation*}%
and hence $\left[ M^{\left( n\right) }-M\right] _{\tau _{k}^{n}}\underset{%
n\rightarrow \infty }{\overset{\mathbb{P}}{\longrightarrow }}0$ for all $%
k\in \mathbb{N}.$

Finally, to prove that $\left[ M^{\left( n\right) }-M\right] _{\infty }%
\overset{\mathbb{P}}{\rightarrow }0$ we assume that it is not true, and then
$\left[ M^{\left( n\right) }-M\right] _{\infty }\overset{\mathbb{P}}{%
\nrightarrow }0$ implies that there exist $\varepsilon >0$ and $\left\{
n_{k}\right\} _{k\in \mathbb{N}}\subset \mathbb{N}$ with
\begin{equation*}
d\left( \left[ M^{\left( n_{k}\right) }-M\right] _{\infty },0\right) \geq
\varepsilon
\end{equation*}%
for all $k\in \mathbb{N},$where $d\left( X,Y\right) :=\inf \left\{
\varepsilon >0:\mathbb{P}\left[ \left\vert X-Y\right\vert >\varepsilon %
\right] \leq \varepsilon \right\} $ is the Ky Fan metric. We shall denote
the subsequence as the original sequence, trying to keep the notation as
simple as possible. Using a diagonal argument, a subsequence $\left\{
n_{i}\right\} _{i\in \mathbb{N}}\subset \mathbb{N}$ can be chosen, with the
property that $d\left( \left[ M^{\left( n_{i}\right) }-M\right] _{\tau
_{k}^{n_{i}}},0\right) <\frac{1}{k}$ for all $i\geq k.$ Since
\begin{equation*}
\lim_{k\rightarrow \infty }\left[ M^{\left( n_{i}\right) }-M\right] _{\tau
_{k}^{n_{i}}}=\left[ M^{\left( n_{i}\right) }-M\right] _{\infty }\quad
\mathbb{P}-a.s.,
\end{equation*}%
we can find some $k\left( n_{i}\right) \geq i$ such that
\begin{equation*}
d\left( \left[ M^{\left( n_{i}\right) }-M\right] _{\tau _{k\left(
n_{i}\right) }^{n_{i}}},\left[ M^{\left( n_{i}\right) }-M\right] _{\infty
}\right) <\frac{1}{k}.
\end{equation*}%
Then, using the estimation
\begin{equation*}
\mathbb{P}\left[ \left\vert \left[ M^{\left( n_{k}\right) }-M\right] _{\tau
_{k\left( n_{k}\right) }^{n_{k}}}-\left[ M^{\left( n_{k}\right) }-M\right]
_{\tau _{k}^{n_{k}}}\right\vert >\varepsilon \right] \leq \mathbb{P}\left[
\left\{ \sup\limits_{t\in \mathbb{R}_{+}}\left\vert M_{t}^{\left(
n_{k}\right) }-M_{t}\right\vert \geq k\right\} \right] ,
\end{equation*}%
it follows that
\begin{equation*}
d\left( \left[ M^{\left( n_{k}\right) }-M\right] _{\tau _{k\left(
n_{k}\right) }^{n_{k}}},\left[ M^{\left( n_{k}\right) }-M\right] _{\tau
_{k}^{n_{k}}}\right) \underset{k\rightarrow \infty }{\longrightarrow }0,
\end{equation*}%
which yields a contradiction with $\varepsilon \leq d\left( \left[ M^{\left(
n_{k}\right) }-M\right] _{\infty },0\right) $. Thus, $\left[ M^{\left(
n\right) }-M\right] _{\infty }\overset{\mathbb{P}}{\rightarrow }0.$
The last part of the this lemma follows immediately from the first statement.
 \hfill $\Box $

Using the Doob's stopping theorem we can conclude that  for $M\in \mathcal{M}_{\infty }$
and an stopping time $\tau $, that $M^{\tau }\in \mathcal{M}_{\infty },$
and therefore it follows as a corollary the following result.

\begin{corollary}
\label{E[|(Mn-M)thau|]->0=>[Mn-M]thau->0_in_P}For $\left\{ M^{\left(
n\right) }\right\} _{n\in \mathbb{N}}\subset \mathcal{M}_{\infty }$, $M\in 
\mathcal{M}_{\infty }$ and $\tau $ any stopping time holds%
\begin{equation*}
M_{\tau }^{\left( n\right) }\overset{L^{1}}{\rightarrow }M_{\tau
}\Longrightarrow \left[ M^{\left( n\right) }-M\right] _{\tau }\overset{%
\mathbb{P}}{\longrightarrow }0.
\end{equation*}
\end{corollary}

\noindent \textit{Proof.}
$\left[ \left( M^{\left( n\right) }\right) ^{\tau }-M^{\tau }\right]
_{\infty }=\left[ M^{\left( n\right) }-M\right] _{\infty }^{\tau }=\left[
M^{\left( n\right) }-M\right] _{\tau }\overset{\mathbb{P}}{\longrightarrow }%
0.$
 \hfill $\Box $

\subsection{Density processes \label{Sect. Density_Processes}}

Given an absolutely continuous probability measure $\mathbb{Q}\ll \mathbb{P}$
in a filtered probability space, where a semimartingale with the weak
predictable representation property is defined, the structure of the density
process has been studied extensively by several authors; see Theorem 14.41
in He, Wang and Yan \cite{HeWanYan} or Theorem III.5.19 in Jacod and
Shiryaev \cite{Jcd&Shry 2003}.

Denote by $D_{t}:=\mathbb{E}\left[ \left. \frac{d\mathbb{Q}}{d\mathbb{P}}%
\right\vert \mathcal{F}_{t}\right] $ the c\`{a}dl\`{a}g version of the
density process. For the increasing sequence of stopping times $\tau
_{n}:=\inf \left\{ t\geq 0:D_{t}<\frac{1}{n}\right\} $ $n\geq 1$ and $\tau
_{0}:=\sup_{n}\tau _{n}$ we have $D_{t}\left( \omega \right) =0$ $\forall
t\geq \tau _{0}\left( \omega \right) $ and $D_{t}\left( \omega \right) >0$ $%
\forall t<\tau _{0}\left( \omega \right) ,$ i.e.%
\begin{equation}
D=D\mathbf{1}_{[\hspace{-0.05cm}[0,\tau _{0}[\hspace{-0.04cm}[},
\label{D=D1[[0,To[[}
\end{equation}%
and the process
\begin{equation}
\frac{1}{D_{s-}}\mathbf{1}_{[\hspace{-0.05cm}[D_{-}\not=0]\hspace{-0.04cm}]}%
\text{ is integrable w.r.t. }D,  \label{1/D_integrable_wrt_D}
\end{equation}%
where we abuse of the notation by setting $[\hspace{-0.05cm}[D_{-}\not=0]%
\hspace{-0.04cm}]:=\left\{ \left( \omega ,t\right) \in \Omega \times \mathbb{%
R}_{+}:D_{t-}\left( \omega \right) \neq 0\right\} .$ Both conditions $\left( %
\ref{D=D1[[0,To[[}\right) $ and $\left( \ref{1/D_integrable_wrt_D}\right) $
are necessary and sufficient in order that a semimartingale to be an \textit{%
exponential semimartigale} \cite[Thm. 9.41]{HeWanYan}, i.e. $D=\mathcal{E}\left( Z\right) $ the Dol\'{e}%
ans-Dade exponential of another semimartingale $Z$. In that case we have
\begin{equation}
\tau _{0}=\inf \left\{ t>0:D_{t-}=0\text{ or }D_{t}=0\right\} =\inf \left\{
t>0:\triangle Z_{t}=-1\right\}.  \label{Tau0=JumpZ=-1}
\end{equation}

It is well known that the L\'{e}vy-processes satisfy the weak property of
predictable representation \cite{HeWanYan}, when the completed natural filtration is
considered. In the following lemma we present the characterization of the
density processes for the case of these processes.

\begin{lemma}
\label{Q<<P =>} Given an absolutely continuous probability measure $\mathbb{Q%
}\ll \mathbb{P}$, there exist coefficients $\theta _{0}\in \mathcal{L}\left(
W\right) \ $and $\theta _{1}\in \mathcal{G}\left( \mu \right) $ such that
\begin{equation}
\frac{d\mathbb{Q}_{t}}{d\mathbb{P}_{t}}=\frac{d\mathbb{Q}_{t}}{d\mathbb{P}%
_{t}}\mathbf{1}_{[\hspace{-0.05cm}[0,\tau _{0}[\hspace{-0.04cm}[}=\mathcal{E}%
\left( Z^{\theta }\right) \left( t\right) ,  \label{Dt=exp(Zt)}
\end{equation}%
where $Z_{t}^{\theta }\in \mathcal{M}_{loc}$ is the local martingale given by%
\begin{equation}
Z_{t}^{\theta }:=\int\limits_{]0,t]}\theta _{0}dW+\int\limits_{]0,t]\times
\mathbb{R}_{0}}\theta _{1}\left( s,x\right) \left( \mu \left( ds,dx\right)
-ds\ \nu \left( dx\right) \right) ,  \label{Def._Ztheta(t)}
\end{equation}%
and $\mathcal{E}$ represents the Doleans-Dade exponential of a
semimartingale. The coefficients $\theta _{0}$ and $\theta _{1}$ are $dt$%
-a.s and $\mu _{\mathbb{P}}^{\mathcal{P}}\left( ds,dx\right) $-a.s. unique
on $[\hspace{-0.05cm}[0,\tau _{0}]\hspace{-0.04cm}]$ and $[\hspace{-0.05cm}%
[0,\tau _{0}]\hspace{-0.04cm}]\times \mathbb{R}_{0}$ respectively for $%
\mathbb{P}$-almost all $\omega $. Furthermore, the coefficients can be
choosen with $\theta _{0}=0$ on $]\hspace{-0.05cm}]\tau _{0},\infty \lbrack
\hspace{-0.04cm}[$ and $\theta _{1}=0$ on $]\hspace{-0.05cm}]\tau
_{0},\infty \lbrack \hspace{-0.04cm}[\times \mathbb{R}$ .
\end{lemma}

\noindent \textit{Proof.} We only address the uniqueness of the coefficients
$\theta _{0}$ and $\theta _{1},$ because the representation follows from $%
\left( \ref{D=D1[[0,To[[}\right) $ and $\left( \ref{1/D_integrable_wrt_D}%
\right) .$ Let assume, that we have two possible vectors $\theta :=\left(
\theta _{0},\theta _{1}\right) $ and $\theta ^{\prime }:=\left( \theta
_{0}^{\prime },\theta _{1}^{\prime }\right) $ satisfying the representation,
i.e.
\begin{equation*}
\begin{array}{rl}
D_{u}\mathbf{1}_{[\hspace{-0.05cm}[0,\tau _{0}[\hspace{-0.04cm}[} & =\int
D_{t-}d\{\int\limits_{]0,t]}\theta _{0}\left( s\right)
dW_{s}+\int\limits_{]0,t]\times \mathbb{R}_{0}}\theta _{1}\left( s,x\right)
\left( \mu \left( ds,dx\right) -ds\ \nu \left( dx\right) \right) \} \\ 
& =\int D_{t-}d\{\int\limits_{]0,t]}\theta _{0}^{\prime }\left( s\right)
dW_{s}+\int\limits_{]0,t]\times \mathbb{R}_{0}}\theta _{1}^{\prime }\left(
s,x\right) \left( \mu \left( ds,dx\right) -ds\ \nu \left( dx\right) \right)
\},%
\end{array}%
\end{equation*}%
and thus%
\begin{eqnarray*}
\triangle D_{t} &=&D_{t-}\triangle \left( \int\limits_{]0,t]\times \mathbb{R%
}_{0}}\theta _{1}\left( s,x\right) \left( \mu \left( ds,dx\right) -ds\ \nu
\left( dx\right) \right) \right)  \\
&=&D_{t-}\triangle \left( \int\limits_{]0,t]\times \mathbb{R}_{0}}\theta
_{1}^{\prime }\left( s,x\right) \left( \mu \left( ds,dx\right) -ds\ \nu
\left( dx\right) \right) \right) .
\end{eqnarray*}% 
Since $D_{t-}>0$ on $[\hspace{-0.05cm}[0,\tau _{0}[\hspace{-0.04cm}[,$
it follows that
\begin{equation*}
\triangle \left( \int\limits_{]0,t]\times \mathbb{R}_{0}}\theta _{1}\left(
s,x\right) \left( \mu \left( ds,dx\right) -ds\ \nu \left( dx\right) \right)
\right) =\triangle \left( \int\limits_{]0,t]\times \mathbb{R}_{0}}\theta
_{1}^{\prime }\left( s,x\right) \left( \mu \left( ds,dx\right) -ds\ \nu
\left( dx\right) \right) \right) .
\end{equation*}
Since two purely discontinuous local martingales with the same jumps are
equal, it follows
\begin{equation*}
\int\limits_{]0,t]\times \mathbb{R}_{0}}\theta _{1}\left( s,x\right) \left(
\mu \left( ds,dx\right) -ds\ \nu \left( dx\right) \right)
=\int\limits_{]0,t]\times \mathbb{R}_{0}}\widehat{\theta }_{1}\left(
s,x\right) \left( \mu \left( ds,dx\right) -ds\ \nu \left( dx\right) \right)
\end{equation*}%
and thus
\begin{equation*}
\int D_{t-}d\{\int\limits_{]0,t]}\theta _{0}\left( s\right) dW_{s}\}=\int
D_{t-}d\{\int\limits_{]0,t]}\theta _{0}^{\prime }\left( s\right) dW_{s}\}.
\end{equation*}%
Then,
\begin{equation*}
0=\left[ \int D_{s-}d\left\{ \int\nolimits_{]0,s]}\left( \theta
_{0}^{\prime }\left( u\right) -\theta _{0}\left( u\right) \right)
dW_{u}\right\} \right] _{t}=\int\limits_{]0,t]}\left( D_{s-}\right)
^{2}\left\{ \theta _{0}^{\prime }\left( s\right) -\theta _{0}\left( s\right)
\right\} ^{2}ds
\end{equation*}%
and thus $\theta _{0}^{\prime }=\theta _{0}\ dt$-$a.s$ on $[\hspace{-0.05cm}%
[0,\tau _{0}]\hspace{-0.04cm}]$ for $\mathbb{P}$-almost all $\omega $.

On the other hand,
\begin{eqnarray*}
0 &=&\left\langle \int \left\{ \theta _{1}^{\prime }\left( s,x\right)
-\theta _{1}\left( s,x\right) \right\} \left( \mu \left( ds,dx\right) -ds\
\nu \left( dx\right) \right) \right\rangle _{t} \\
&=&\int\limits_{]0,t]\times \mathbb{R}_{0}}\left\{ \theta _{1}^{\prime
}\left( s,x\right) -\theta _{1}\left( s,x\right) \right\} ^{2}\nu \left(
dx\right) ds,
\end{eqnarray*}%
implies that $\theta _{1}\left( s,x\right) =\theta _{1}^{\prime }\left(
s,x\right) \quad \mu _{\mathbb{P}}^{\mathcal{P}}\left( ds,dx\right) $-a.s.
on $[\hspace{-0.05cm}[0,\tau _{0}]\hspace{-0.04cm}]\times \mathbb{R}_{0}$
for $\mathbb{P}$-almost all $\omega $. \hfill $\Box $

For $\mathbb{Q}\ll \mathbb{P}$ the function $\theta _{1}\left( \omega
,t,x\right) $ described in Lemma \ref{Q<<P =>} determines the density of the
predictable projection $\mu _{\mathbb{Q}}^{\mathcal{P}}\left( dt,dx\right) $
with respect to $\mu _{\mathbb{P}}^{\mathcal{P}}\left( dt,dx\right) $ (see
He,Wang and Yan \cite{HeWanYan} or Jacod and Shiryaev \cite{Jcd&Shry 2003}).
More precisely, for\ $B\in \left( \mathcal{B}\left( \mathbb{R}_{+}\right)
\otimes \mathcal{B}\left( \mathbb{R}_{0}\right) \right) $ we have
\begin{equation}
\mu _{\mathbb{Q}}^{\mathcal{P}}\left( \omega ,B\right) =\int_{B}\left(
1+\theta _{1}\left( \omega ,t,x\right) \right) \mu _{\mathbb{P}}^{\mathcal{P}%
}\left( dt,dx\right) .  \label{Q<<P=>_miu_wrt_Q}
\end{equation}

In what follows we restrict ourself to the time interval $\left[ 0,T\right]
, $ for some $T>0$ fixed, and we take $\mathcal{F}=\mathcal{F}_{T}.$ The
corresponding classes of density processes associated to $\mathcal{Q}_{\ll }(%
\mathbb{P})$ and $\mathcal{Q}_{\approx }\left( \mathbb{P}\right) $ are
denoted by $\mathcal{D}_{\ll }\left( \mathbb{P}\right) $ and $\mathcal{D}%
_{\approx }\left( \mathbb{P}\right) $, respectively. For instance, in the
former case
\begin{equation}
\mathcal{D}_{\ll }\left( \mathbb{P}\right) :=\left\{ D=\left\{ D_{t}\right\}
_{t\in \left[ 0,T\right] }:\exists \mathbb{Q}\in \mathcal{Q}_{\ll }\left(
\mathbb{P}\right) \text{ with }D_{t}=\left. \frac{d\mathbb{Q}}{d\mathbb{P}}%
\right\vert _{\mathcal{F}_{t}}\right\} ,  \label{Def._D<<}
\end{equation}%
and the processes in this set are of the form
\begin{equation}
\begin{array}{rl}
D_{t}= & \exp \left\{ \int\limits_{]0,t]}\theta
_{0}dW+\int\limits_{]0,t]\times \mathbb{R}_{0}}\theta _{1}\left( s,x\right)
\left( \mu \left( ds,dx\right) -\nu \left( dx\right) ds\right) -\frac{1}{2}%
\int\limits_{]0,t]}\left( \theta _{0}\right) ^{2}ds\right\} \times \\
& \times \exp \left\{ \int\limits_{]0,t]\times \mathbb{R}_{0}}\left\{ \ln
\left( 1+\theta _{1}\left( s,x\right) \right) -\theta _{1}\left( s,x\right)
\right\} \mu \left( ds,dx\right) \right\}%
\end{array}
\label{D(t) explicita}
\end{equation}%
for $\theta _{0}\in \mathcal{L}\left( W\right) $ and $\theta _{1}\in
\mathcal{G}\left( \mu \right) $.

The set $\mathcal{D}_{\ll }\left( \mathbb{P}\right) $ is characterized as
follow.

\begin{corollary}
\label{D<<_<=>} $D$ belongs to $\mathcal{D}_{\ll }\left( \mathbb{P}\right) $
if and only if there are $\theta _{0}\in \mathcal{L}\left( W\right) $ and $%
\theta _{1}\in \mathcal{G}\left( \mu \right) $ with $\theta _{1}\geq -1$
such that $D_{t}=\mathcal{E}\left( Z^{\theta }\right) \left( t\right) \
\mathbb{P}$-a.s. $\forall t\in \left[ 0,T\right] $ and $\mathbb{E}_{\mathbb{P%
}}\left[ \mathcal{E}\left( Z^{\theta }\right) \left( t\right) \right] =1\
\forall t\geq 0$, where $Z^{\theta }\left( t\right) $ is defined by $\left( %
\ref{Def._Ztheta(t)}\right) .$
\end{corollary}

\noindent \textit{Proof.} The necessity follows from Lemma \ref{Q<<P =>}.
Conversely, let $\theta _{0}\in \mathcal{L}\left( W\right) $ and $\theta
_{1}\in \mathcal{G}\left( \mu \right) $ be arbitrarily chosen. Since $%
D_{t}=\int D_{s-}dZ_{s}^{\theta }\in \mathcal{M}_{loc}$ is a nonnegative
local martingale, it is a supermartingale, with constant expectation from
our assumptions. Therefore, it is a martingale, and hence the density
process of an absolutely continuous probability measure. \hfill $\Box$

Since density processes are essentially uniformly integrable martingales,
using Lemma \ref{E[|Mn-M|]->0=>[Mn-M](oo)->0_in_P} and Corollary %
\ref{E[|(Mn-M)thau|]->0=>[Mn-M]thau->0_in_P}   the
following proposition follows immediately.

\begin{proposicion}
\label{E[|Dn-D|]->0 => [Dn-D](T)->0_in_P} Let $\left\{ \mathbb{Q}^{\left(
n\right) }\right\} _{n\in \mathbb{N}}$ be a sequence in $\mathcal{Q}_{\ll }(%
\mathbb{P})$, with $D_{T}^{\left( n\right) }:=\left. \frac{d\mathbb{Q}%
^{\left( n\right) }}{d\mathbb{P}}\right\vert _{\mathcal{F}_{T}}$ converging
to $D_{T}:=\left. \frac{d\mathbb{Q}}{d\mathbb{P}}\right\vert _{\mathcal{F}%
_{T}}$ in $L^{1}\left( \mathbb{P}\right) $. For the corresponding density
processes $D_{t}^{\left( n\right) }:=\mathbb{E}_{\mathbb{P}}\left[
D_{T}^{\left( n\right) }\left\vert \mathcal{F}_{t}\right. \right] $ and $%
D_{t}:=\mathbb{E}_{\mathbb{P}}\left[ D_{T}\left\vert \mathcal{F}_{t}\right. %
\right] $, for $t\in \left[ 0,T\right] $, we have%
\begin{equation*}
\left[ D^{\left( n\right) }-D\right] _{T}\overset{\mathbb{P}}{\rightarrow }0.
\end{equation*}
\end{proposicion} 

\section{Penalty functions for densities\label{Sect Penalty Function for
densities}}

Now, we shall introduce a family of penalty functions for the density
processes described in Section \ref{Sect. Density_Processes}, for the
absolutely continuous measures $\mathbb{Q}\in \mathcal{Q}_{\ll }\left(
\mathbb{P}\right) $.

Let $h:\mathbb{R}_{+}\mathbb{\rightarrow R}_{+}$ and $h_{0},$\thinspace $%
h_{1}:\ \mathbb{R\rightarrow R}_{+}$ be convex functions with $0=h\left(
0\right) =h_{0}\left( 0\right) =h_{1}\left( 0\right) $. Define the penalty
function, with $\tau_0$ as in (\ref{Tau0=JumpZ=-1}), by
\begin{equation}
\begin{array}{rl}
\vartheta \left( \mathbb{Q}\right) := & \mathbb{E}_{\mathbb{Q}}\left[
\int\limits_{0}^{T\wedge \tau _{0}}h\left( h_{0}\left( \theta _{0}\left(
t\right) \right) +\int\nolimits_{\mathbb{R}_{0}}\delta \left( t,x\right)
h_{1}\left( \theta _{1}\left( t,x\right) \right) \nu \left( dx\right)
\right) dt\right] \mathbf{1}_{\mathcal{Q}_{\ll }}\left( \mathbb{Q}\right) \\
& +\infty \times \mathbf{1}_{\mathcal{Q}_{cont}\setminus \mathcal{Q}_{\ll
}}\left( \mathbb{Q}\right) ,%
\end{array}
\label{Def._penalty_theta}
\end{equation}%
 where $\theta _{0},$ $\theta _{1}$ are the processes associated to $\mathbb{Q%
}$ from Lemma \ref{Q<<P =>} and $\delta \left( t,x\right) :\mathbb{R}%
_{+}\times \mathbb{R}_{0}\rightarrow \mathbb{R}_{+}$ is an arbitrary fixed
nonnegative function $\delta \left( t,x\right) \in \mathcal{G}\left( \mu
\right) $. Since $\theta _{0}\equiv 0$ on $[\hspace{-0.05cm}[\tau
_{0},\infty \lbrack \hspace{-0.04cm}[$ and $\theta _{1}\equiv 0$ on $[%
\hspace{-0.05cm}[\tau _{0},\infty \lbrack \hspace{-0.04cm}[\times \mathbb{R}%
_{0}$ we have from the conditions imposed to $h,h_{0},$ and $h_{1}$%
\begin{equation}
\begin{array}{rl}
\vartheta \left( \mathbb{Q}\right) = & \mathbb{E}_{\mathbb{Q}}\left[
\int\limits_{0}^{T}h\left( h_{0}\left( \theta _{0}\left( t\right) \right)
+\int\nolimits_{\mathbb{R}_{0}}\delta \left( t,x\right) h_{1}\left( \theta
_{1}\left( t,x\right) \right) \nu \left( dx\right) \right) dt\right] \mathbf{%
1}_{\mathcal{Q}_{\ll }}\left( \mathbb{Q}\right) \\
& +\infty \times \mathbf{1}_{\mathcal{Q}_{cont}\setminus \mathcal{Q}_{\ll
}}\left( \mathbb{Q}\right) .%
\end{array}
\label{Def._penalty_theta_(2)}
\end{equation}%
Further, define the convex measure of risk
\begin{equation}
\rho \left( X\right) :=\sup_{\mathbb{Q\in }\mathcal{Q}_{\ll }(\mathbb{P}%
)}\left\{ \mathbb{E}_{\mathbb{Q}}\left[ -X\right] -\vartheta \left( \mathbb{Q%
}\right) \right\} .  \label{rho def.}
\end{equation}%
Notice that $\rho $ is a normalized and sensitive measure of risk. For each
class of probability measures introduced so far, the subclass of those
measures with a finite penalization is considered. We will denote by $%
\mathcal{Q}^{\vartheta },$ $\mathcal{Q}_{\ll }^{\vartheta }(\mathbb{P})$ and
$\mathcal{Q}_{\approx }^{\vartheta }(\mathbb{P})$ the respective subclasses,
i.e.
\begin{equation}
\mathcal{Q}^{\vartheta }:=\left\{ \mathbb{Q}\in \mathcal{Q}:\vartheta \left(
\mathbb{Q}\right) <\infty \right\} ,\ \mathcal{Q}_{\ll }^{\vartheta }(%
\mathbb{P}):=\mathcal{Q}^{\vartheta }\cap \mathcal{Q}_{\ll }(\mathbb{P})%
\text{ and }\mathcal{Q}_{\approx }^{\vartheta }(\mathbb{P}):=\mathcal{Q}%
^{\vartheta }\cap \mathcal{Q}_{\approx }(\mathbb{P}).  \label{Def._Qdelta(P)}
\end{equation}%
Notice that $\mathcal{Q}_{\approx }^{\vartheta
}(\mathbb{P})\neq \varnothing .$

The next theorem establishes the minimality on $\mathcal{Q}_{\ll }\left(
\mathbb{P}\right) $ of the penalty function introduced above for the risk
measure $\rho $ . Its proof is based on the sufficient conditions given in
Theorem \ref{static minimal penalty funct. in Q(<<) <=>}.

\begin{theorem}
\label{theta=minimal penalty function} The penalty function $\vartheta $
defined in $\left( \ref{Def._penalty_theta}\right) $ is equal to    the minimal penalty function of the convex
risk measure $\rho $, given by $\left( \ref{rho def.}\right) $, on $\mathcal{Q}%
_{\ll }\left( \mathbb{P}\right) $, i.e.%
\begin{equation*}
\vartheta \mathbf{1}_{\mathcal{Q}_{\ll }\left( \mathbb{P}\right) }=\psi
_{\rho }^{\ast }\mathbf{1}_{\mathcal{Q}_{\ll }\left( \mathbb{P}\right) }.
\end{equation*}
\end{theorem}

\noindent \textit{Proof:} From Lemma \ref{static minimal penalty funct. in
Q(<<) <=>} $\left( b\right) $, we need to show that the penalization $%
\vartheta $ is proper, convex and that the corresponding identification,
defined as $\Theta \left( Z\right) :=\vartheta \left( \mathbb{Q}\right) $ if
$Z\mathbb{\in }\delta \left( \mathcal{Q}_{\ll }\left( \mathbb{P}\right)
\right) :=\left\{ Z\in L^{1}\left( \mathbb{P}\right) :Z=d\mathbb{Q}/d\mathbb{%
P}\text{ with }\mathbb{Q}\in \mathcal{Q}_{\ll }\left( \mathbb{P}\right)
\right\} $ and $\Theta \left( Z\right) :=\infty $ on $L^{1}\setminus \delta
\left( \mathcal{Q}_{\ll }\left( \mathbb{P}\right) \right) $, is lower
semicontinuous with respect to the strong topology.

First, observe that the function $\vartheta $ is proper, since $\vartheta
\left( \mathbb{P}\right) =0$. To verify the convexity of $\vartheta $,
choose $\mathbb{Q}$, $\widetilde{\mathbb{Q}}\in \mathcal{Q}_{\ll
}^{\vartheta }$ and define $\mathbb{Q}^{\lambda }:=\lambda \mathbb{Q}+\left(
1-\lambda \right) \widetilde{\mathbb{Q}}$, for $\lambda \in \left[ 0,1\right]
$. Notice that the corresponding density process can be written as $%
D^{\lambda }:=\dfrac{d\mathbb{Q}^{\lambda }}{d\mathbb{P}}=\lambda D+\left(
1-\lambda \right) \widetilde{D}$ $\mathbb{P}$-a.s. .

Now, from Lemma \ref{Q<<P =>}, let $\left( \theta _{0},\theta _{1}\right) $
and $(\widetilde{\theta }_{0},\widetilde{\theta }_{1})$ be the processes
associated to $\mathbb{Q}$ and $\widetilde{\mathbb{Q}}$, respectively, and
observe that from%
\begin{equation*}
D_{t}=1+\int\limits_{\left[ 0,t\right] }D_{s-}\theta _{0}\left( s\right)
dW_{s}+\int\limits_{\left[ 0,t\right] \times \mathbb{R}_{0}}D_{s-}\theta
_{1}\left( s,x\right) d\left( \mu \left( ds,dx\right) -ds\nu \left(
dx\right) \right) )
\end{equation*}%
and the corresponding expression for $\widetilde{D}$ we have for $\tau
_{n}^{\lambda }:=\inf \left\{ t\geq 0:D_{t}^{\lambda }\leq \frac{1}{n}%
\right\} $
\begin{equation*}
\int\limits_{0}^{t\wedge \tau _{n}^{\lambda }}\left( D_{s-}^{\lambda
}\right) ^{-1}dD_{s}^{\lambda }=\int\limits_{0}^{t\wedge \tau _{n}^{\lambda
}}\tfrac{\lambda D_{s-}\theta _{0}\left( s\right) +\left( 1-\lambda \right)
\widetilde{D}_{s-}\widetilde{\theta }_{0}\left( s\right) }{\left( \lambda
D_{s-}+\left( 1-\lambda \right) \widetilde{D}_{s-}\right) }%
dW_{s}+\int\limits_{\left[ 0,t\wedge \tau _{n}^{\lambda }\right] \times
\mathbb{R}_{0}}\tfrac{\lambda D_{s-}\theta _{1}\left( s,x\right) +\left(
1-\lambda \right) \widetilde{D}_{s-}\widetilde{\theta }_{1}\left( s,x\right)
}{\left( \lambda D_{s-}+\left( 1-\lambda \right) \widetilde{D}_{s-}\right) }%
d\left( \mu -\mu _{\mathbb{P}}^{\mathcal{P}}\right) .
\end{equation*}%
The weak predictable representation property of the local martingale $%
\int\nolimits_{0}^{t\wedge \tau _{n}^{\lambda }}\left( D_{s-}^{\lambda
}\right) ^{-1}dD_{s}^{\lambda }$, yield on the other hand
\begin{equation*}
\int\limits_{0}^{t\wedge \tau _{n}^{\lambda }}\left( D_{s-}^{\lambda
}\right) ^{-1}dD_{s}^{\lambda }=\int\limits_{0}^{t\wedge \tau _{n}^{\lambda
}}\theta _{0}^{\lambda }\left( s\right) dW_{s}+\int\limits_{\left[ 0,t\wedge
\tau _{n}^{\lambda }\right] \times \mathbb{R}_{0}}\theta _{1}^{\lambda
}\left( s,x\right) d\left( \mu -\mu _{\mathbb{P}}^{\mathcal{P}}\right) ,
\end{equation*}%
where identification
\begin{equation*}
\theta _{0}^{\lambda }\left( s\right) =\frac{\lambda D_{s-}\theta _{0}\left(
s\right) +\left( 1-\lambda \right) \widetilde{D}_{s-}\widetilde{\theta }%
_{0}\left( s\right) }{\left( \lambda D_{s-}+\left( 1-\lambda \right)
\widetilde{D}_{s-}\right) },
\end{equation*}%
and
\begin{equation*}
\theta _{1}^{\lambda }\left( s,x\right) =\frac{\lambda D_{s-}\theta
_{1}\left( s,x\right) +\left( 1-\lambda \right) \widetilde{D}_{s-}\widetilde{%
\theta }_{1}\left( s,x\right) }{\left( \lambda D_{s-}+\left( 1-\lambda
\right) \widetilde{D}_{s-}\right) }.
\end{equation*}%
\ is possible thanks to the uniqueness of the representation in Lemma \ref%
{Q<<P =>}. The convexity follows now from the convexity of $h,h_{0}$%
\thinspace and $h_{1}$, using the fact that any convex function is
continuous in the interior of its domain. More specifically,
\begin{equation*}
\begin{array}{rl}
\vartheta \left( \mathbb{Q}^{\lambda }\right) \leq & \mathbb{E}_{\mathbb{Q}%
^{\lambda }}\left[ \int\limits_{\left[ 0,T\right] }\tfrac{\lambda D_{s}}{%
\left( \lambda D_{s}+\left( 1-\lambda \right) \widetilde{D}_{s}\right) }%
h\left( h_{0}\left( \theta _{0}\left( s\right) \right) +\int\limits_{\mathbb{%
R}_{0}}\delta \left( s,x\right) h_{1}\left( \theta _{1}\left( s,x\right)
\right) \nu \left( dx\right) \right) ds\right] \\
& +\mathbb{E}_{\mathbb{Q}^{\lambda }}\left[ \int\limits_{\left[ 0,T\right] }%
\tfrac{\left( 1-\lambda \right) \widetilde{D}_{s}}{\left( \lambda
D_{s}+\left( 1-\lambda \right) \widetilde{D}_{s}\right) }h\left( h_{0}\left(
\widetilde{\theta }_{0}\left( s\right) \right) +\int\limits_{\mathbb{R}%
_{0}}\delta \left( s,x\right) h_{1}(\widetilde{\theta }_{1}\left( s,x\right)
)\nu \left( dx\right) \right) ds\right] \\
= & \int\limits_{\left[ 0,T\right] }\int\limits_{\Omega }\dfrac{\lambda D_{s}%
}{\left( \lambda D_{s}+\left( 1-\lambda \right) \widetilde{D}_{s}\right) }%
h\left( h_{0}\left( \theta _{0}\left( s\right) \right) +\int\limits_{\mathbb{%
R}_{0}}\delta \left( s,x\right) h_{1}\left( \theta _{1}\left( s,x\right)
\right) \nu \left( dx\right) \right) \\
& \ \ \ \ \ \ \ \ \times \left( \lambda D_{s}+\left( 1-\lambda \right)
\widetilde{D}_{s}\right) \mathbf{1}_{\left\{ \lambda D_{s}+\left( 1-\lambda
\right) \widetilde{D}_{s}>0\right\} }d\mathbb{P}ds \\
& +\int\limits_{\left[ 0,T\right] }\int\limits_{\Omega }\dfrac{\left(
1-\lambda \right) \widetilde{D}_{s}}{\left( \lambda D_{s}+\left( 1-\lambda
\right) \widetilde{D}_{s}\right) }h\left( h_{0}\left( \widetilde{\theta }%
_{0}\left( s\right) \right) +\int\limits_{\mathbb{R}_{0}}\delta \left(
s,x\right) h_{1}(\widetilde{\theta }_{1}\left( s,x\right) )\nu \left(
dx\right) \right) \\
& \ \ \ \ \ \ \ \ \times \left( \lambda D_{s}+\left( 1-\lambda \right)
\widetilde{D}_{s}\right) \mathbf{1}_{\left\{ \lambda D_{s}+\left( 1-\lambda
\right) \widetilde{D}_{s}>0\right\} }d\mathbb{P}ds \\
= & \lambda \vartheta \left( \mathbb{Q}\right) +\left( 1-\lambda \right)
\vartheta \left( \widetilde{\mathbb{Q}}\right) ,%
\end{array}%
\end{equation*}%
where we used that $\left\{ \int\nolimits_{\mathbb{R}_{0}}\delta \left(
t,x\right) h_{1}\left( \theta _{1}\left( t,x\right) \right) \nu \left(
dx\right) \right\} _{t\in \mathbb{R}_{+}}$ and $\left\{ \int\nolimits_{%
\mathbb{R}_{0}}\delta \left( t,x\right) h_{1}(\widetilde{\theta }_{1}\left(
t,x\right) )\nu \left( dx\right) \right\} _{t\in \mathbb{R}_{+}}$ are
predictable processes.

It remains to prove the lower semicontinuity of $\Theta $. As pointed out
earlier, it is enough to consider a sequence of densities $Z^{\left(
n\right) }:=\frac{d\mathbb{Q}^{\left( n\right) }}{d\mathbb{P}}\in \delta
\left( \mathcal{Q}_{\ll }\left( \mathbb{P}\right) \right) $ converging in $%
L^{1}\left( \mathbb{P}\right) $ to $Z:=\frac{d\mathbb{Q}}{d\mathbb{P}}$.
Denote the corresponding density processes by $D^{\left( n\right) }$ and $D$%
, respectively. In Proposition \ref{E[|Dn-D|]->0 => [Dn-D](T)->0_in_P} it was
verified the convergence in probability to zero of the quadratic variation
process
\begin{eqnarray*}
\left[ D^{\left( n\right) }-D\right] _{T} &=&\int\limits_{0}^{T}\left\{
D_{s-}^{\left( n\right) }\theta _{0}^{\left( n\right) }\left( s\right)
-D_{s-}\theta _{0}\left( s\right) \right\} ^{2}ds \\
&&+\int\limits_{\left[ 0,T\right] \times \mathbb{R}_{0}}\left\{
D_{s-}^{\left( n\right) }\theta _{1}^{\left( n\right) }\left( s,x\right)
-D_{s-}\theta _{1}\left( s,x\right) \right\} ^{2}\mu \left( ds,dx\right) .
\end{eqnarray*}%
This implies that
\begin{equation}
\left.
\begin{array}{cc}
& \int\nolimits_{0}^{T}\left\{ D_{s-}^{\left( n\right) }\theta _{0}^{\left(
n\right) }\left( s\right) -D_{s-}\theta _{0}\left( s\right) \right\} ^{2}ds%
\overset{\mathbb{P}}{\rightarrow }0, \\
\text{and } &  \\
& \int\limits_{\left[ 0,T\right] \times \mathbb{R}_{0}}\left\{
D_{s-}^{\left( n\right) }\theta _{1}^{\left( n\right) }\left( s,x\right)
-D_{s-}\theta _{1}\left( s,x\right) \right\} ^{2}\mu \left( ds,dx\right)
\overset{\mathbb{P}}{\rightarrow }0.%
\end{array}%
\right\}  \label{[]=>*}
\end{equation}%
Then, for an arbitrary but fixed subsequence, there exists a sub-subsequence
such that $\mathbb{P}$-a.s.
\begin{equation*}
\left\{ D_{s-}^{\left( n\right) }\theta _{0}^{\left( n\right) }\left(
s\right) -D_{s-}\theta _{0}\left( s\right) \right\} ^{2}\overset{L^{1}\left(
\lambda \right) }{\longrightarrow }0
\end{equation*}%
and
\begin{equation*}
\left\{ D_{s-}^{\left( n\right) }\theta _{1}^{\left( n\right) }\left(
s,x\right) -D_{s-}\theta _{1}\left( s,x\right) \right\} ^{2}\overset{%
L^{1}\left( \mu \right) }{\longrightarrow }0,
\end{equation*}%
where for simplicity we have denoted the sub-subsequence as the original
sequence. Now, we claim that for the former sub-subsequence it also holds
that
\begin{equation}
\left\{
\begin{array}{c}
D_{s-}^{\left( n\right) }\theta _{0}^{\left( n\right) }\left( s\right)
\overset{\lambda \times \mathbb{P}\text{-a.s.}}{\longrightarrow }%
D_{s-}\theta _{0}\left( s\right) , \\
\smallskip \  \\
D_{s-}^{\left( n\right) }\theta _{1}^{\left( n\right) }\left( s,x\right)
\overset{\mu \times \mathbb{P}\text{-a.s.}}{\longrightarrow }D_{s-}\theta
_{1}\left( s,x\right) .%
\end{array}%
\right.  \label{[]=>*.1}
\end{equation}

We present first the arguments for the proof of the second assertion in $%
\left( \ref{[]=>*.1}\right) $. Assuming the opposite, there exists $C\in
\mathcal{B}\left( \left[ 0,T\right] \right) \otimes \mathcal{B}\left(
\mathbb{R}_{0}\right) \otimes \mathcal{F}_{T}$, with $\mu \times \mathbb{P}%
\left[ C\right] >0$, and such that for each $\left( s,x,\omega \right) \in C$
\begin{equation*}
\lim_{n\rightarrow \infty }\left\{ D_{s-}^{\left( n\right) }\theta
_{1}^{\left( n\right) }\left( s,x\right) -D_{s-}\theta _{1}\left( s,x\right)
\right\} ^{2}=c\neq 0,
\end{equation*}%
or the limit does not exist.

Let $C\left( \omega \right) :=\left\{ \left( t,x\right) \in \left[ 0,T\right]
\times \mathbb{R}_{0}:\left( t,x,\omega \right) \in C\right\} $ be the $%
\omega $-section of $C$. Observe that $B:=\left\{ \omega \in \Omega :\mu %
\left[ C\left( \omega \right) \right] >0\right\} $ has positive probability:
$\mathbb{P}\left[ B\right] >0.$

From $\left( \ref{[]=>*}\right) $, any arbitrary but fixed subsequence has a
sub-subsequence converging $\mathbb{P}$-a.s..  Denoting such a
sub-subsequence simply by $n$, we can fix $\omega \in B$ with%
\begin{eqnarray*}
&&\int\nolimits_{C\left( \omega \right) }\left\{ D_{s-}^{\left( n\right)
}\theta _{1}^{\left( n\right) }\left( s,x\right) -D_{s-}\theta _{1}\left(
s,x\right) \right\} ^{2}d\mu \left( s,x\right) \\
&\leq &\int\nolimits_{\left[ 0,T\right] \times \mathbb{R}_{0}}\left\{
D_{s-}^{\left( n\right) }\theta _{1}^{\left( n\right) }\left( s,x\right)
-D_{s-}\theta _{1}\left( s,x\right) \right\} ^{2}d\mu \left( s,x\right)
\underset{n\rightarrow \infty }{\longrightarrow }0,
\end{eqnarray*}%
and hence $\left\{ D_{s-}^{\left( n\right) }\theta _{1}^{\left( n\right)
}\left( s,x\right) -D_{s-}\theta _{1}\left( s,x\right) \right\} ^{2}$
converges in $\mu $-measure to $0$ on $C\left( \omega \right) .$ Again, for
any subsequence there is a sub-subsequence converging $\mu $-a.s. to $0$.
Furthermore, for an arbitrary but fixed $\left( s,x\right) \in C\left(
\omega \right) $, when the limit does not exist
\begin{equation*}
\begin{array}{clc}
a & :=\underset{n\rightarrow \infty }{\lim \inf }\left\{ D_{s-}^{\left(
n\right) }\theta _{1}^{\left( n\right) }\left( s,x\right) -D_{s-}\theta
_{1}\left( s,x\right) \right\} ^{2} &  \\
& \neq \underset{n\rightarrow \infty }{\lim \sup }\left\{ D_{s-}^{\left(
n\right) }\theta _{1}^{\left( n\right) }\left( s,x\right) -D_{s-}\theta
_{1}\left( s,x\right) \right\} ^{2} & =:b,%
\end{array}%
\end{equation*}%
and we can choose converging subsequences $n\left( i\right) $ and $n\left(
j\right) $ with
\begin{eqnarray*}
\underset{i\rightarrow \infty }{\lim }\left\{ D_{s-}^{n\left( i\right)
}\theta _{1}^{n\left( i\right) }\left( s,x\right) -D_{s-}\theta _{1}\left(
s,x\right) \right\} ^{2} &=&a \\
\underset{j\rightarrow \infty }{\lim }\left\{ D_{s-}^{n\left( j\right)
}\theta _{1}^{n\left( j\right) }\left( s,x\right) -D_{s-}\theta _{1}\left(
s,x\right) \right\} ^{2} &=&b.
\end{eqnarray*}%
From the above argument, there are sub-subsequences $n\left( i\left(
k\right) \right) $ and $n\left( j\left( k\right) \right) $ such that
\begin{eqnarray*}
a &=&\underset{k\rightarrow \infty }{\lim }\left\{ D_{s-}^{n\left( i\left(
k\right) \right) }\theta _{1}^{n\left( i\left( k\right) \right) }\left(
s,x\right) -D_{s-}\theta _{1}\left( s,x\right) \right\} ^{2}=0 \\
b &=&\underset{k\rightarrow \infty }{\lim }\left\{ D_{s-}^{n\left( j\left(
k\right) \right) }\theta _{1}^{n\left( j\left( k\right) \right) }\left(
s,x\right) -D_{s-}\theta _{1}\left( s,x\right) \right\} ^{2}=0,
\end{eqnarray*}%
which is clearly a contradiction.

For the case when
\begin{equation*}
\underset{n\rightarrow \infty }{\lim }\left\{ D_{s-}^{\left( n\right)
}\theta _{1}^{\left( n\right) }\left( s,x\right) -D_{s-}\theta _{1}\left(
s,x\right) \right\} ^{2}=c\neq 0,
\end{equation*}%
the same argument can be used, and get a subsequence converging to $0$,
having a contradiction again. Therefore, the second part of our claim in $%
\left( \ref{[]=>*.1}\right) $ holds.

Since $D_{s-}^{\left( n\right) }\theta _{1}^{\left( n\right) }\left(
s,x\right) ,\ D_{s-}\theta _{1}\left( s,x\right) \in \mathcal{G}\left( \mu
\right) $, we have, in particular, that $D_{s-}^{\left( n\right) }\theta
_{1}^{\left( n\right) }\left( s,x\right) \in \widetilde{\mathcal{P}}$ and $%
D_{s-}\theta _{1}\left( s,x\right) \in \widetilde{\mathcal{P}}$ and hence $%
C\in \widetilde{\mathcal{P}}$. From the definition of the predictable
projection it follows that
\begin{eqnarray*}
0 &=&\mu \times \mathbb{P}\left[ C\right] \mathbb{=}\int\limits_{\Omega
}\int\limits_{\left[ 0,T\right] \times \mathbb{R}_{0}}\mathbf{1}_{C}\left(
s,\omega \right) d\mu d\mathbb{P=}\int\limits_{\Omega }\int\limits_{\left[
0,T\right] \times \mathbb{R}_{0}}\mathbf{1}_{C}\left( s,\omega \right) d\mu
_{\mathbb{P}}^{\mathcal{P}}d\mathbb{P} \\
&=&\int\limits_{\Omega }\int\limits_{\mathbb{R}_{0}}\int\limits_{\left[ 0,T%
\right] }\mathbf{1}_{C}\left( s,\omega \right) dsd\nu d\mathbb{P=}\lambda
\times \nu \times \mathbb{P}\left[ C\right] ,
\end{eqnarray*}%
and thus
\begin{equation*}
D_{s-}^{\left( n\right) }\theta _{1}^{\left( n\right) }\left( s,x\right)
\overset{\lambda \times \nu \times \mathbb{P}\text{-a.s.}}{\longrightarrow }%
D_{s-}\theta _{1}\left( s,x\right) .
\end{equation*}

Since $\int\limits_{\Omega \times \left[ 0,T\right] }\left\vert
D_{t-}^{\left( n\right) }-D_{t-}\right\vert d\mathbb{P}\times dt\mathbb{=}%
\int\limits_{\Omega \times \left[ 0,T\right] }\left\vert D_{t}^{\left(
n\right) }-D_{t}\right\vert d\mathbb{P}\times dt\longrightarrow 0$, we have
that\linebreak\ $\left\{ D_{t-}^{\left( n\right) }\right\} _{t\in \left[ 0,T%
\right] }$ $\overset{L^{1}\left( \lambda \times \mathbb{P}\right) }{%
\longrightarrow }\left\{ D_{t-}\right\} _{t\in \left[ 0,T\right] }$ and $%
\left\{ D_{t}^{\left( n\right) }\right\} _{t\in \left[ 0,T\right] }$ $%
\overset{L^{1}\left( \lambda \times \mathbb{P}\right) }{\longrightarrow }%
\left\{ D_{t}\right\} _{t\in \left[ 0,T\right] }.$ Then, for an arbitrary
but fixed subsequence $\left\{ n_{k}\right\} _{k\in \mathbb{N}}\subset
\mathbb{N}$, there is a sub-subsequence $\left\{ n_{k_{i}}\right\} _{i\in
\mathbb{N}}\subset \mathbb{N}$ such that
\begin{equation*}
\begin{array}{ccc}
D_{t-}^{\left( n_{k_{i}}\right) }\theta _{1}^{\left( n_{k_{i}}\right)
}\left( t,x\right) & \overset{\lambda \times \nu \times \mathbb{P}\text{-a.s.%
}}{\longrightarrow } & D_{t-}\theta _{1}\left( t,x\right) , \\
D_{t-}^{\left( n_{k_{i}}\right) } & \overset{\lambda \times \mathbb{P}\text{%
-a.s.}}{\longrightarrow } & D_{t-}, \\
D_{t}^{\left( n_{k_{i}}\right) } & \overset{\lambda \times \mathbb{P}\text{%
-a.s.}}{\longrightarrow } & D_{t}.%
\end{array}%
\end{equation*}%
Furthermore, $\mathbb{Q}\ll \mathbb{P}$ implies that $\lambda \times \nu
\times \mathbb{Q}\ll \lambda \times \nu \times \mathbb{P}$, and then
\begin{equation*}
\begin{array}{ccc}
D_{t-}^{\left( n_{k_{i}}\right) }\theta _{1}^{\left( n_{k_{i}}\right)
}\left( t,x\right) & \overset{\lambda \times \nu \times \mathbb{Q}\text{-a.s.%
}}{\longrightarrow } & D_{t-}\theta _{1}\left( t,x\right) , \\
D_{t-}^{\left( n_{k_{i}}\right) } & \overset{\lambda \times \nu \times
\mathbb{Q}\text{-a.s.}}{\longrightarrow } & D_{t-},%
\end{array}%
\end{equation*}%
and
\begin{equation}
D_{t}^{\left( n_{k_{i}}\right) }\overset{\lambda \times \nu \times \mathbb{Q}%
\text{-a.s.}}{\longrightarrow }D_{t}.  \label{[]=>*.2}
\end{equation}%
Finally, noting that $\inf D_{t}>0$ $\mathbb{Q}$-a.s.
\begin{equation}
\theta _{1}^{\left( n_{k_{i}}\right) }\left( t,x\right) \overset{\lambda
\times \nu \times \mathbb{Q}\text{-a.s.}}{\longrightarrow }\theta _{1}\left(
t,x\right) .  \label{[]=>*.3}
\end{equation}

The first assertion in $\left( \ref{[]=>*.1}\right) $ can be proved using
essentially the same kind of ideas used above for the proof of the second
part, concluding that for an arbitrary but fixed subsequence $\left\{
n_{k}\right\} _{k\in \mathbb{N}}\subset \mathbb{N}$, there is a
sub-subsequence $\left\{ n_{k_{i}}\right\} _{i\in \mathbb{N}}\subset \mathbb{%
N}$ such that
\begin{equation}
\left\{ D_{t}^{\left( n_{k_{i}}\right) }\right\} _{t\in \left[ 0,T\right] }%
\overset{\lambda \times \mathbb{Q}\text{-a.s.}}{\longrightarrow }\left\{
D_{t}\right\} _{t\in \left[ 0,T\right] }  \label{[]=>*.4}
\end{equation}%
and
\begin{equation}
\left\{ \theta _{0}^{\left( n_{k_{i}}\right) }\left( t\right) \right\}
_{t\in \left[ 0,T\right] }\overset{\lambda \times \mathbb{Q}\text{-a.s.}}{%
\longrightarrow }\left\{ \theta _{0}\left( t\right) \right\} _{t\in \left[
0,T\right] }.  \label{[]=>*.5}
\end{equation}

We are now ready to finish the proof of the theorem, observing that
\begin{eqnarray*}
&&\underset{n\rightarrow \infty }{\lim \inf }\vartheta \left( \mathbb{Q}%
^{\left( n\right) }\right) \\
&=&\underset{n\rightarrow \infty }{\lim \inf }\int\limits_{\Omega \times %
\left[ 0,T\right] }\left\{ h\left( h_{0}\left( \theta _{0}^{\left( n\right)
}\left( t\right) \right) +\int\nolimits_{\mathbb{R}_{0}}\delta \left(
t,x\right) h_{1}\left( \theta _{1}^{\left( n\right) }\left( t,x\right)
\right) \nu \left( dx\right) \right) \right\} \dfrac{D_{t}^{\left( n\right) }%
}{D_{t}}d\left( \lambda \times \mathbb{Q}\right) .
\end{eqnarray*}%
Let $\left\{ n_{k}\right\} _{k\in \mathbb{N}}\subset \mathbb{N}$ be a
subsequence for which the limit inferior is realized. Using $\left( \ref%
{[]=>*.2}\right) ,\left( \ref{[]=>*.3}\right) ,\ $\linebreak $\left( \ref%
{[]=>*.4}\right) ,$ and $\left( \ref{[]=>*.5}\right) $ we can pass to a
sub-subsequence $\left\{ n_{k_{i}}\right\} _{i\in \mathbb{N}}\subset \mathbb{%
N}$ and, from the continuity of $h,\ h_{0}$\thinspace and\thinspace $h_{1}$,
it follows
\begin{eqnarray*}
&&\underset{n\rightarrow \infty }{\lim \inf }\ \vartheta \left( \mathbb{Q}%
^{\left( n\right) }\right) \\
&\geq &\int\limits_{\Omega \times \left[ 0,T\right] }\underset{i\rightarrow
\infty }{\lim \inf }\left( \left\{ h\left( h_{0}\left( \theta _{0}^{\left(
n_{k_{i}}\right) }\left( t\right) \right) +\int\limits_{\mathbb{R}%
_{0}}\delta \left( t,x\right) h_{1}\left( \theta _{1}^{\left(
n_{k_{i}}\right) }\left( t,x\right) \right) \nu \left( dx\right) \right)
\right\} \tfrac{D_{t}^{\left( n_{k_{i}}\right) }}{D_{t}}\right) d\left(
\lambda \times \mathbb{Q}\right) \\
&\geq &\int\limits_{\Omega \times \left[ 0,T\right] }h\left( h_{0}\left(
\theta _{0}\left( t\right) \right) +\int\nolimits_{\mathbb{R}%
_{0}}h_{1}\left( \theta _{1}\left( t,x\right) \right) \nu \left( dx\right)
\right) d\left( \lambda \times \mathbb{Q}\right) \\
&=&\vartheta \left( \mathbb{Q}\right) .
\end{eqnarray*}%
\hfill $\Box $

\eject

\end{document}